\documentclass[conference,twoside,10pt]{ieeeconf}
\usepackage{cite}
\usepackage{xcolor}
\usepackage{amsmath,amssymb,amsfonts}
\usepackage{comment}
\usepackage{graphicx}
\usepackage{mathrsfs}
\usepackage{empheq}
\newtheorem{theorem}{Theorem}
\newtheorem{remark}{Remark}
\newtheorem{definition}{Definition}
\newtheorem{lem}{Lemma}

\newenvironment{definitionn}{\begin{definition}}{\hfill $\square$ \end{definition}}
\newenvironment{lemma}{\begin{lem}}{\hfill $\square$ \end{lem}}
\usepackage{hyperref}

\allowdisplaybreaks
\usepackage{textcomp}
\begin{document}
\title{\bf $H^2$ Stabilization of the $2$-D and $3$-D\\ Heat Equation via Modal Decomposition}
\author{Mohamed Amine Ouchdiri, Mohamed Camil Belhadjoudja, Mohamed Maghenem,\\ Saad Benjelloun, Adnane Saoud
\thanks{M.~A. Ouchdiri and A.~Saoud are with College of
Computing, University Mohammed VI Polytechnic, Benguerir, Morocco (e-mail: amine.ouchdiri,adnane.saoud@um6p.ma).}
\thanks{M. C. Belhadjoudja and M.~Maghenem are with University Grenoble Alpes, CNRS, Grenoble-INP, GIPSA-lab, F-38000, Grenoble, France (e-mail: mohamed.belhadjoudja,mohamed.maghenem@gipsa-lab.fr).}
\thanks{S.~Benjelloun is with De Vinci Higher Education, De Vinci Research Center, Paris, France (e-mail: saad.benjelloun@devinci.fr).}}
\maketitle
\pagenumbering{gobble}

\begin{abstract}
 Modal decomposition has been used 
 recently to achieve boundary 
 stabilization of linear parabolic equations in the $H^1$ norm, in two and three dimensions. 
 In one dimension, $H^1$ exponential stability is well established and known to imply boundedness and asymptotic convergence of the state's max norm to zero.  However, this implication does not necessarily hold in higher dimensions. In this paper, we prove that the 
 modal-decomposition-based full-state feedback controller in \cite{Munteanu2017IJC} guarantees, not only $H^1$ exponential stability, but also
 $H^2$ exponential stability, which implies boundedness and asymptotic convergence of the state's max norm to zero. 
 Our approach is based on rewriting the state's Laplacian, required in the expression of the state's $H^2$ norm, as a linear combination of the state and its time derivative. Thereby, the $L^2$ norm of the state's time derivative is the one to be analyzed in terms of boundedness and asymptotic convergence to zero.
\end{abstract}

\section{Introduction}\label{sec1}

Parabolic partial differential equations (PDEs) model tissue-scale signaling and pattern formation in developmental biology \cite{moh2,moh}, heat conduction, and atmospheric pollutant dispersion. One-dimensional models usually suffice in some simplified settings, while most physical processes require higher-dimensional domains.

Despite the prevalence of 
high-dimensional parabolic PDEs, 
few methods are proposed to stabilize a given equilibrium using boundary control.  Infinite-dimensional optimal control is considered in  \cite{LasieckaTriggiani2000}, relying on operator Riccati equations that may be numerically intractable. Infinite-dimensional backstepping  on symmetric domains or parallelepipeds is considered in \cite{VazquezKrstic2016,LiuXie2020,Meurer2013}, on some directionally-convex domains in \cite{higher_1}, and on more general domains but for specific boundary control configurations and  dynamic controllers in \cite{higher_2}. 
Related to our framework,  
approximation-based methods, via modal decomposition, are considered in  \cite{Balas1978,KatzFridman2020,KatzFridman2021,LhachemiPrieur2022}, where the stability guarantees for a finite-dimensional linear approximation are extended to the original PDE. Most existing results based on modal decomposition consider one-dimensional PDEs only. 
The first extension to higher dimensions is proposed in \cite{Barbu2013TAC}, assuming the independence of some of  the eigenfunctions' normal derivatives. This assumption was later relaxed in \cite{Munteanu2017IJC} and \cite{lhachemi2025automatica} using  full-state- and output-feedback controllers, respectively. 

The strongest stability guarantee obtained using a modal-decomposition framework is in the sense of the $H^1$ norm. In the one-dimensional setting, $H^1$ exponential stability is known to imply boundedness and asymptotic convergence of the state's max norm to zero. However, this implication is not systematic in higher dimensions. 
Motivated by this fact, we propose to prove $H^2$ exponential stability in the case of full-state feedback, and we show that such a property does imply boundedness and convergence of the state's max norm.

Boundary stabilization in the $H^2$ norm is established in \cite{lhachemi2024stabilization,LhachemiPrieur2026H2}, for a one-dimensional reaction--diffusion equation. 
The analysis therein uses a functional, that is equivalent to the squared $H^2$ norm. Since the differentiation of such functional produces an infinite series that is not convergent in general, instead, a truncation of such a functional is used as a Lyapunov function candidate. This approach is difficult to extend to higher dimensions, especially when the actuation is of Dirichlet-type, which we consider here. This suggests taking a different route. Indeed, observing that the $H^1$ stability is already established in \cite{lhachemi2025automatica}, we conclude that it remains to study the $L^2$ norm of the state's Laplacian. However, the heat equation gives the state's Laplacian as a linear combination of the state and its time derivative. Moreover, the $L^2$ norm of the state is upper-bounded by its $H^1$ norm. Hence, it is enough to study the $L^2$ norm of the state's time derivative.  We analyze the latter by differentiating, with respect to time, the equations governing the equation's modes.

We introduce preliminary results in Section~\ref{sec2}, the control system in Section~\ref{sec3}, our main result in Section~\ref{sec4}, and numerical simulations in Section~\ref{sec6}.

\textit{Notations.} For $n\in \mathbb{N}^*:=\{1,2,...\}$ and $v := (v_1,v_2,...,v_n)\in \mathbb{R}^n$, we let $|v|:=\sqrt{v_1^2+v_2^2+...+v_n^2}$ denote the usual Euclidean norm on $\mathbb{R}^n$. Moreover, for $w:=(w_1,w_2,...,w_n)$, we let $\langle v,w\rangle_n := v_1w_1+v_2w_2+...+v_nw_n$ denote the inner product associated to $|\cdot|$. Furthermore, given a matrix $P\in \mathbb{R}^{n\times n}$, we let $\|P\|:=\sup_{\{a\in \mathbb{R}^{n}, \ |a|=1\}}\{|Pa|\}$ be the induced matrix norm. Moreover, given a bounded open set $\Omega\subset \mathbb{R}^n$, we denote by $L^2(\Omega)$ the space of (classes of equivalence of) functions $f:\Omega\to \mathbb{R}$ such that $\|f\|_{L^2(\Omega)} := \sqrt{\int_{\Omega}f(x)^2dx}<+\infty$ with $x=(x_1,x_2,...,x_n)\in \Omega$. The inner product on $L^2(\Omega)$ is given by $\langle f, g \rangle_{L^2(\Omega)} := \int_{\Omega}f(x)g(x)dx$ for all $f,g\in L^2(\Omega)$. Furthermore, we let $H^1(\Omega)$ denote the space of functions $f\in L^2(\Omega)$ such that $|\nabla f|\in L^2(\Omega)$, where $\nabla f := \left(f_{x_1}, f_{x_2},...,f_{x_n}\right)$ is the gradient of $f$, and $f_{x_i}$, for $i\in \{1,2,...,n\}$, is the partial derivative of $f$ with respect to $x_i$. The norm on $H^1(\Omega)$ is given by $\|f\|_{H^1(\Omega)} := \sqrt{\|f\|_{L^2(\Omega)}^2+\||\nabla f|\|_{L^2(\Omega)}^2}$. Additionally, we let $H_0^1(\Omega):=\{f\in H^1(\Omega):f|_{\partial\Omega}=0\}$, where $\partial \Omega$ denotes the boundary of $\Omega$, and $f|_{\partial \Omega}$ denotes the restriction of $f$ to $\partial \Omega$. Similarly, we let $H^2(\Omega)$ be the space of functions $f\in H^1(\Omega)$ such that $\{f_{x_ix_j}\}_{i,j=1}^{n}\subset L^2(\Omega)$, where $f_{x_ix_j}$, for $i,j\in \{1,2,...,n\}$, denotes the partial derivative of $f_{x_i}$ with respect to $x_j$. Furthermore, we denote by $\Delta f := f_{x_1x_1}+f_{x_2x_2}+...+f_{x_nx_n}$ the Laplacian of $f$. The norm on $H^2(\Omega)$ is denoted $\|f\|_{H^2(\Omega)} := \sqrt{\|f\|_{H^1(\Omega)}^2+\sum_{i,j=1}^{n}\|f_{x_ix_j}\|_{L^2(\Omega)}^2}$. We denote by $H^{\frac{3}{2}}(\partial\Omega)$ the image of $H^{2}(\Omega)$ under the
Dirichlet trace operator $\omega_{0}(f):=f|_{\partial\Omega}$, equipped with the
quotient norm
\begin{equation*}
  \|g\|_{H^{\frac{3}{2}}(\partial\Omega)}
  := \inf\Big\{\,\|f\|_{H^{2}(\Omega)} \;:\; f\in H^{2}(\Omega),\
  f|_{\partial\Omega}=g \,\Big\}.
\end{equation*}
For a continuous function $f:\bar{\Omega} \to \mathbb{R}$ where $\bar{\Omega}$ is the closure of $\Omega$, we let $\|f\|_{L^{\infty}(\Omega)} := \max_{x\in \bar{\Omega}}\{|f(x)|\}$. Finally, we denote by $\mathcal{C}^k([0,+\infty);X)$ the space of $k$-times continuously differentiable functions $f:[0,+\infty)\to X$, $t\mapsto f(t)$, with $\partial_t f \in X$ the derivative of $f$ with respect to $t$, and $\partial_t^2f := \partial_t \partial_tf$.

\section{Preliminaries}\label{sec2}
In the sequel, we let $\Omega \subset \mathbb{R}^n$, with $n \in \{1,2,3\}$, be a bounded domain with closure $\bar{\Omega}$. The boundary of $\Omega$ is denoted by $\partial \Omega$, and is assumed to be of class $\mathcal{C}^3$ \cite[p.~710]{Evans2010}. In this paper, we study the problem of $H^2$ stabilization of a linear heat equation posed on $\Omega$, achieved via boundary control and modal decomposition. The control system and our objective will be formulated precisely in Section \ref{sec3}. At this stage, we introduce some operators that will play a central role throughout the paper and summarize their main properties.  

\subsection{The Heat Operator} 

Let $\lambda \geq 0$ and $\mathcal{A}:D(\mathcal{A})\to L^2(\Omega)$ be defined by
\begin{align}
\mathcal{A} \varphi := \Delta \varphi + \lambda \varphi \quad \forall \varphi \in D(\mathcal{A}),  \label{eq4}
\end{align}
whose domain is given by
\begin{align}
D(\mathcal{A}) := \{\varphi \in H^2(\Omega) \ : \ \varphi = 0 \ \text{on $\partial \Omega$}\}.\label{eq5}
\end{align}
Since $\Omega$ is a bounded domain with a $\mathcal{C}^3$ boundary, one has the following result \cite[Thm.~1, p.~355]{Evans2010}. 
\begin{lemma}
The operator $\mathcal{A}$ has a countable set of eigenvalues $\{\mu_{i}\}_{i\in \mathbb{N}^*}\subset \mathbb{R}$, ordered as $\mu_1\geq \mu_2 \geq \cdots \geq \mu_i \geq \cdots$, with $\mu_i\to -\infty$ as $i\to +\infty$, and associated $L^2(\Omega)$-orthonormal eigenfunctions $\{\varphi_i\}_{i\in\mathbb{N}^*}$ forming a Riesz basis of $L^2(\Omega)$. 
\end{lemma}

\subsection{The Boundary Homogenizing Operator}

Let $N\in\mathbb{N}^*$ denote the number of nonnegative eigenvalues of $\mathcal{A}$, i.e., 
\begin{align*}
\mu_1\geq \mu_2 \geq \cdots \geq \mu_N \geq 0 > \mu_{N+1} \geq \cdots
\end{align*}

Given $\gamma >0$ and $f\in H^{\frac{3}{2}}(\partial \Omega)$, we consider the equation
\begin{equation*}
\Sigma_H : \left\lbrace
\begin{aligned}
&-\bar{\mathcal{A}}D+2\sum_{i=1}^{N}\mu_i \langle D,\varphi_i\rangle_{L^2(\Omega)}\varphi_i+\gamma D = 0 \quad \text{on $\Omega$}, \\
&D = f \quad \text{on $\partial \Omega$}, 
\end{aligned}
\right.
\end{equation*}
where $\bar{\mathcal{A}} : H^2(\Omega)\to L^2(\Omega)$ is defined by $\bar{\mathcal{A}}f := \Delta f + \lambda f$. Note that, since $D$ equals $f$ on $\partial \Omega$, subtracting $D$ from any function that also equals $f$ on $\partial \Omega$ gives a function that is zero on $\partial \Omega$. This operator will be essential later on to homogenize the considered control system, a step that is standard when employing modal decomposition. 

Given the $\mathcal{C}^3$ regularity of $\partial \Omega$, one has the following key result, whose proof is in Appendix \ref{p_lem_D}. 

\begin{lemma}\label{lem_D}
The following properties hold. 
\begin{itemize}
    \item[(i)] For a constant $\gamma>0$ large enough and for every $f\in H^{\frac{3}{2}}(\partial \Omega)$, there exists a unique function $D\in H^2(\Omega)$ that satisfies $\Sigma_H$.
    \item[(ii)] For $\gamma$ large enough, define the operator 
\begin{align}
D_{\gamma}:H^{\frac{3}{2}}(\partial \Omega)\to H^2(\Omega) \label{D_g}
\end{align}
which maps any $f\in H^{\frac{3}{2}}(\partial \Omega)$ to the solution $D_{\gamma}(f) := D$ of $\Sigma_H$. Then, $D_{\gamma}:H^{\frac{3}{2}}(\partial \Omega)\to H^2(\Omega)$ is linear and continuous, and
\begin{align*}
\|D_{\gamma}(f)\|_{H^2(\Omega)} \leq C_{D,\gamma} \|f\|_{H^{\frac{3}{2}}(\partial \Omega)} \quad \forall f \in H^{\frac{3}{2}}(\partial \Omega),
\end{align*}
where 
$
C_{D,\gamma} := \sup_{\substack{f \in H^{\frac{3}{2}}(\partial \Omega) \\ f \neq 0}} \frac{\|D_{\gamma}(f)\|_{H^2(\Omega)}}{\|f\|_{H^{\frac{3}{2}}(\partial \Omega)}} < +\infty.
$
\item[(iii)] Let $v \in \mathcal{C}^1([0,+\infty); H^{\frac{3}{2}}(\partial \Omega))$. Then, $D_{\gamma}(v(t))\in \mathcal{C}^1([0,+\infty); H^2(\Omega))$ and
$\partial_t D_{\gamma}(v(t)) = D_{\gamma}\left( \partial_t v(t) \right).$
\end{itemize}
\end{lemma}

\subsection{The Normal Trace Operator} 
Let the operator $T_{\mathbf{n}} : H^3(\Omega)\to H^{\frac{3}{2}}(\partial \Omega)$, defined by
\begin{align*}
T_{\mathbf{n}}(\varphi) := \langle \nabla\varphi, \mathbf{n}\rangle_n \quad \forall \varphi\in H^3(\Omega), 
\end{align*}
where $\mathbf{n}$ denotes the unit outward normal to $\partial \Omega$. 

Note that $T_{\mathbf{n}}$ can be applied to the eigenfunctions of $\mathcal{A}$. Indeed, since $\mathcal{A}\varphi_i=\mu_i\varphi_i\in H^1(\Omega)$, one has $\varphi_i\in H^3(\Omega)$ by \cite[Thm.~2.5.1.1, p.~128]{Grisvard}, and $T_{\mathbf{n}}(\varphi_i)\in H^{\frac{3}{2}}(\partial\Omega)$ by \cite[Thm.~1.5.1.2, pp.~37--38]{Grisvard}.

\section{Problem Formulation}\label{sec3}

We consider the heat equation
\begin{empheq}[left=\Sigma:\ \left\{,right=\right.]{align}
&\partial_t u = \Delta u + \lambda u \quad &&\text{on $\Omega \times (0,+\infty)$}, \label{eq1_heat}\\
&u(x,t) = v(x,t) \quad &&\text{on $\partial \Omega \times [0,+\infty)$}, \label{eq2_heat}\\
&u(x,0)= u_o(x) \quad &&\text{on $\bar{\Omega}$}, \label{eq3_heat}
\end{empheq}
where $u:\bar{\Omega}\times [0,+\infty)\to \mathbb{R}$ is the state, $x\in \bar{\Omega}$ is the space variable, $t\geq 0$ is the time variable, $v(x,t)\in \mathbb{R}$ is the control variable, $u_o\in H^2(\Omega)$ is the initial condition, and $\lambda \geq 0$ is the reaction coefficient.

The solutions to $\Sigma$ are understood in the strong sense, recalled below.

\begin{definitionn}\label{def1}
Given $u_o\in H^2(\Omega)$, a strong solution to $\Sigma$ is a function
$u\in \mathcal{C}([0,+\infty);H^2(\Omega))\cap \mathcal{C}^1([0,+\infty);L^2(\Omega))$
that satisfies \eqref{eq1_heat}, \eqref{eq2_heat}, and \eqref{eq3_heat} almost everywhere.
\end{definitionn}

Under $v=0$, $t \mapsto \|u(\cdot,t)\|_{L^2(\Omega)}$ is unbounded 
if $\lambda$ is large enough. This has motivated, in \cite{Munteanu2017IJC,lhachemi2025automatica}, the design of $v$ to achieve exponential stability of the origin $\{u=0\}$ for $\Sigma$ in $H^1$. However, in $2$-D and $3$-D, $H^1$ exponential stability does not imply pointwise convergence of the state to zero nor its boundedness. Hence, our objective here is to prove the latter two properties by guaranteeing $H^2$ exponential stability, thanks to the following key Lemma \cite[Thms. 1.5.1.2 and 2.4.2.5]{Grisvard}.

\begin{lemma}\label{lem2}
Let $\Omega\subset\mathbb{R}^n$, with $n\in\{1,2,3\}$, be a bounded and connected domain with $\mathcal{C}^3$ boundary. There exist constants $C_{S_1},C_{S_2}>0$, depending on $\Omega$ but independent of $f$, such that, for all $f\in H^2(\Omega)$,
\begin{align}
&\|f\|_{H^2(\Omega)} \leq C_{S_1}\big(\|f\|_{L^2(\Omega)}+\|\Delta f\|_{L^2(\Omega)}\notag\\
&\hspace{1.6cm}+\|f|_{\partial\Omega}\|_{H^{\frac{3}{2}}(\partial\Omega)}\big), \label{eq11}\\
&\|f\|_{L^\infty(\Omega)} \leq C_{S_2}\|f\|_{H^2(\Omega)}. \label{eq12}
\end{align}
\end{lemma} 
\begin{remark}
Inequality \eqref{eq12} is false in dimensions higher than $3$. This is the reason why we require that $\Omega \subset \mathbb{R}^n$ with $n\leq 3$ in order to establish boundedness and pointwise convergence of the state to zero. 
\end{remark}

\subsection{The Controller in \cite{Munteanu2017IJC}}\label{sec3b}
As standard in modal decomposition, the first step is to homogenize $\Sigma$, i.e., to perform a change of variables $w = \mathcal{T}(u)$, for some transformation $\mathcal{T}$ to be defined below, so that $w$ vanishes on $\partial \Omega$. Among the various possible choices, the appropriate definition of $w$ is as follows\footnote{The main contribution of \cite{Munteanu2017IJC} was to identify this specific change of variables (cf. \eqref{eq16}--\eqref{eq15}), thereby enabling the subsequent steps of the modal decomposition procedure. We refer the reader to the aforementioned reference for a justification of this specific transformation over alternative choices. In essence, other transformations considered in the literature yield, at some stage of the modal decomposition procedure, a matrix that is singular whenever the normal derivatives of the eigenfunctions of $\mathcal{A}$, corresponding to the eigenvalues $\mu_1,...,\mu_N$, are linearly dependent \cite{Barbu2013TAC}.}:
\begin{itemize}
\item We select $N$ constants $0 < \gamma_1 < \gamma_2 < \cdots < \gamma_N$, all sufficiently large and distinct from the eigenvalues $\mu_1, \ldots, \mu_N$, so that each operator $D_{\gamma_i}$ is well-defined. 

\item Then, we decompose the control input $v$ as
\begin{align}
v(x,t) = \sum_{i=1}^{N} v_i(x,t), \label{eq14}
\end{align}
where each $v_i : \partial\Omega \times [0,+\infty) \to \mathbb{R}$ is a function to be designed. 

\item For each $i \in \{1,\ldots,N\}$, we apply the lifting operator $D_{\gamma_i}$ to $v_i$, and let
\begin{align}
\xi_i(x,t) := D_{\gamma_i}(v_i(\cdot,t)). \label{eq16}
\end{align}

\item Finally, we define
\begin{align}
w = \mathcal{T}(u) := u - \sum_{i=1}^{N}\xi_i. \label{eq15}
\end{align}
\end{itemize}

By construction, each $\xi_i$ satisfies $\xi_i = v_i$ on $\partial\Omega$. Summing over $i$ and using \eqref{eq14}, we obtain $\sum_{i=1}^{N}\xi_i = v$ on $\partial\Omega$. Hence, $w$ vanishes on $\partial \Omega$.

One verifies that $w$ solves the PDE 
\begin{equation*}
\Sigma_w : \left\lbrace 
\begin{aligned}
&\partial_t w = \mathcal{A}w -\sum_{i=1}^{N}\partial_t \xi_i + \sum_{i=1}^{N}\gamma_i \xi_i \\
&\hspace{0.4cm}+ 2\sum_{k,i=1}^{N}\mu_k \langle \xi_i,\varphi_k\rangle_{L^2(\Omega)} \varphi_k \ \ \text{on } \Omega \times (0,+\infty), \\
&w(x,t) = 0 \ \ \text{on } \partial \Omega \times [0,+\infty), \\
&w(x,0)=u_o(x)-\sum_{i=1}^{N}\xi_i(x,0) \ \ \text{on } \bar{\Omega}. 
\end{aligned}
\right.
\end{equation*}

The next step is to project $u$ and $w$ on the eigenfunctions of $\mathcal{A}$. That is, we let, for each $i\in\mathbb{N}^*$,
\begin{align*}
u_i := \langle u,\varphi_i\rangle_{L^2(\Omega)}, \quad w_i := \langle w,\varphi_i\rangle_{L^2(\Omega)}.
\end{align*}
Note that $u_i$ and $w_i$ can be linked, through \eqref{eq15}, as follows
\begin{align}
u_i = w_i + \sum_{k=1}^{N}\langle \xi_k,\varphi_i\rangle_{L^2(\Omega)} \quad \forall\, i\in \{1,\ldots,N\}. \label{eq17}
\end{align}

Using \eqref{eq17}, we can show that the \textit{vector of unstable modes} $U(t) := [u_1(t),\ldots,u_N(t)]^\top \in \mathbb{R}^N$ satisfies
\begin{align}
&\frac{d}{dt}U(t) = A_oU(t) - 2\sum_{i=1}^{N}A_oM_{\gamma_i}\begin{bmatrix}
    \langle v_i,T_{\mathbf{n}}(\varphi_1) \rangle_{L^2(\partial \Omega)} \\
    \vdots \\
    \langle v_i,T_{\mathbf{n}}(\varphi_N) \rangle_{L^2(\partial \Omega)}
\end{bmatrix} \nonumber \\
&\quad -\sum_{i=1}^{N}\left(-A_o+\gamma_i I_{N}\right)M_{\gamma_i}\begin{bmatrix}
    \langle v_i,T_{\mathbf{n}}(\varphi_1) \rangle_{L^2(\partial \Omega)} \\
    \vdots \\
    \langle v_i,T_{\mathbf{n}}(\varphi_N) \rangle_{L^2(\partial \Omega)}
\end{bmatrix}, \label{eq18}
\end{align}
where
\begin{align}
A_o &:= \mathrm{diag}(\mu_1,\ldots,\mu_N) \in \mathbb{R}^{N\times N}, \label{eq20} \\
M_{\gamma_i} &:= \mathrm{diag}\!\left(\frac{1}{\gamma_i+\mu_1},\ldots,\frac{1}{\gamma_i+\mu_N}\right) \in \mathbb{R}^{N\times N}. \label{eq21}
\end{align}

Finally, we design $\{v_i\}_{i=1}^N$ to stabilize the origin $\{U=0\}$ for \eqref{eq18}. To this end, for each $i\in\{1,\ldots,N\}$, we set
\begin{align}
v_i(x,t) := \langle M_{\gamma_i}AU(t),\,\mathcal{L}(x)\rangle_N \quad x\in\partial\Omega,\; t\geq 0, \label{eq22}
\end{align}
where
\begin{itemize}
    \item $\mathcal{L}(x) := \big[T_{\mathbf{n}}(\varphi_1)(x),\ldots,T_{\mathbf{n}}(\varphi_N)(x)\big]^\top$,
    \item $A\in\mathbb{R}^{N\times N}$ is a matrix defined as follows: let $B\in\mathbb{R}^{N\times N}$ be the Gram matrix of the normal traces, whose $(i,j)$-entry is
$
[B]_{ij} := \langle T_{\mathbf{n}}(\varphi_i),\,T_{\mathbf{n}}(\varphi_j)\rangle_{L^2(\partial\Omega)},$
and for each $i\in\{1,\ldots,N\}$, define
\begin{align}
B_i := M_{\gamma_i}\,B\,M_{\gamma_i}. \label{eq23}
\end{align}
It is shown in \cite{Munteanu2017IJC} that the matrix $\sum_{i=1}^{N}B_i$ is invertible. We then define
$
A := \left(\sum_{i=1}^{N}B_i\right)^{-1}. 
$
\end{itemize}

Substituting \eqref{eq22} into \eqref{eq18}, we obtain, after simplification,
\begin{align*}
\frac{d}{dt}U(t) = \left(-\sum_{i=1}^{N}\gamma_i B_i A\right)U(t).
\end{align*}
It is shown that if the parameters $\gamma_1,\ldots,\gamma_N$ are chosen sufficiently large, then the matrix $-\sum_{i=1}^{N}\gamma_i B_i A$ is Hurwitz \cite{Munteanu2017IJC}. Hence, there exist $C_1\geq 1$ and $\sigma_o > 0$ such that
\begin{align}
|U(t)| \leq C_1\,|U(0)|\,e^{-\sigma_o t} \quad \forall\, t\geq 0. \label{eq28}
\end{align}

The proof that \eqref{eq28} implies $L^2$ exponential stability of the origin for $\Sigma$ by analyzing the variations of a suitable Lyapunov functional is standard. This analysis is extended in \cite{lhachemi2025automatica} to conclude $H^1$ exponential stability. Precisely, we have the following result.

\begin{lemma}[{\cite{Munteanu2017IJC,lhachemi2025automatica}}]\label{lem3}
Consider the system $\Sigma$, and assume that the operator $\mathcal{A}$ in \eqref{eq4}-\eqref{eq5} has $N\in \mathbb{N}^*$ nonnegative eigenvalues $\mu_1\geq \cdots \geq \mu_N\geq 0$. Moreover, let the boundary input $v$ be given by \eqref{eq14} and \eqref{eq22}, where the parameters $0<\gamma_1<\cdots<\gamma_N$ are chosen distinct from $\mu_1,\ldots, \mu_N$, and sufficiently large such that
\begin{itemize}
    \item For each $i\in \{1,\ldots,N\}$, the equation $\Sigma_H$ with $\gamma$ replaced by $\gamma_i$, and $f\in H^{\frac{3}{2}}(\partial \Omega)$, admits a unique solution $D=D_{\gamma_i}(f) \in H^2(\Omega)$.
    \item The matrix $-\sum_{i=1}^{N}\gamma_i B_iA$, where  $B_i$, $i\in \{1,\ldots,N\}$,  in \eqref{eq23} and $A :=(\sum_{i=1}^{N}B_i)^{-1}$, is Hurwitz.
\end{itemize}
Then, there exist two constants $\Gamma \geq 1$ and $\sigma>0$ such that, for any initial condition $u_o \in H^2(\Omega)$ satisfying the compatibility condition $u_o|_{\partial \Omega} = v(\cdot,0)$, we have 
\begin{align}
\|u(\cdot,t)\|_{H^1(\Omega)}\leq \Gamma \|u_o\|_{H^1(\Omega)}e^{-\sigma t} \quad \forall t\geq 0. \label{eq29}
\end{align}
\end{lemma}

\section{$H^2$ Exponential Stability}\label{sec4}

We complete Lemma \ref{lem3} into the following main result. 

\begin{theorem}\label{thm_main}
Let the hypotheses of Lemma \ref{lem3} hold. Then,
 there exist constants
$\Gamma_1 \geq 1$ and $\sigma^*>0$, independent of the initial condition $u_o$, such that the solutions to $\Sigma$ satisfy
\begin{align}
    \|u(\cdot,t)\|_{H^2(\Omega)} 
    &\leq \Gamma_1\,\|u_o\|_{H^2(\Omega)}\,
    e^{-\sigma^* t} \quad \forall t\geq 0. \label{eq32}
\end{align}

In particular, there exists a constant $\Gamma_2 \geq 1$, independent of $u_o$, such that
\begin{align}
    \|u(\cdot,t)\|_{L^\infty(\Omega)} 
    &\leq \Gamma_2\,\|u_o\|_{H^2(\Omega)}\,
    e^{-\sigma^* t} \quad \forall t\geq 0. \label{eq33}
\end{align}

\hfill $\square$
\end{theorem}

\subsection{Proof of Theorem \ref{thm_main}}
By Lemma \ref{lem2}, we need to upperbound $\|u(\cdot,t)\|_{L^2(\Omega)}$, $\|\Delta u(\cdot,t)\|_{L^2(\Omega)}$, and $\|u(\cdot,t)|_{\partial\Omega}\|_{H^{\frac{3}{2}}(\partial\Omega)}$. The first term decays exponentially to zero according to Lemma~\ref{lem3}. The third term is just the norm of the control input $v(\cdot,t)$, satisfying
\begin{align}
\|v(\cdot,t)\|_{H^{\frac{3}{2}}(\partial\Omega)} \leq C_{\mathcal{T}}\sum_{i=1}^{N}\|M_{\gamma_i}A\|\,|U(t)|, \label{eq32_bis}
\end{align}
where $C_{\mathcal{T}} := \sum_{j=1}^{N}\|T_{\mathbf{n}}(\varphi_j)\|_{H^{\frac{3}{2}}(\partial\Omega)}$, so that it decays exponentially to zero at the rate $\sigma_o$ by \eqref{eq28}. Hence, it remains to establish the exponential decay of $t \mapsto \|\Delta u(\cdot,t)\|_{L^2(\Omega)}$. To this end, we first note that, from $\Sigma$, one has
\begin{align}
\Delta u = \partial_t u - \lambda u, \label{eq30}
\end{align}
and, hence, by the triangle inequality 
\begin{align}
\hspace{-0.3cm} \|\Delta u(\cdot,t)\|_{L^2(\Omega)} \leq \|\partial_t u(\cdot,t)\|_{L^2(\Omega)} + \lambda\,\|u(\cdot,t)\|_{L^2(\Omega)}. \label{last_eq}
\end{align}
Since $\|u(\cdot,t)\|_{L^2(\Omega)} \leq \|u(\cdot,t)\|_{H^1(\Omega)}$, the second term in the right-hand side of \eqref{last_eq} decays exponentially to zero by Lemma~\ref{lem3}. Therefore, the proof of $H^2$ exponential stability reduces to establishing the exponential decay to zero of $t \mapsto \|\partial_t u(\cdot,t)\|_{L^2(\Omega)}$.
To do so, we differentiate both sides of \eqref{eq15} with respect to time, to obtain
\begin{align}
\partial_t u = \partial_t w + \sum_{i=1}^{N}\partial_t \xi_i. \label{eq31}
\end{align}
Then, we carefully upper-bound the $L^2$ norm of each term on the right-hand side of \eqref{eq31}. Specifically, we prove the following key lemma.
\begin{lemma}\label{lem4}
Under the hypotheses of Lemma~\ref{lem3}, there exist constants $C_\xi,C_w,\sigma^* > 0$, independent of $u_o$, such that, for all $i \in \{1,\ldots,N\}$ and all $t \geq 0$,
\begin{align}
\|\partial_t\xi_i(\cdot,t)\|_{H^2(\Omega)} &\leq C_\xi\,\|u_o\|_{H^2(\Omega)}\,e^{-\sigma_o t}, \label{eq34} \\
\|\xi_i(\cdot,t)\|_{H^2(\Omega)} &\leq C_\xi\,\|u_o\|_{H^2(\Omega)}\,e^{-\sigma_o t}, \label{eq35} \\
\|\partial_t w(\cdot,t)\|_{L^2(\Omega)} &\leq C_w\,\|u_o\|_{H^2(\Omega)}\,e^{-\sigma^* t}. \label{eq38} 
\end{align}
\end{lemma}

Once Lemma \ref{lem4} is proven, we can conclude the proof of Theorem \ref{thm_main}. Indeed, from this lemma and \eqref{eq31}, setting $C_{u_t} := C_w + N\,C_\xi$, we obtain
\begin{align}
\|\partial_t u(\cdot,t)\|_{L^2(\Omega)} \leq C_{u_t}\,\|u_o\|_{H^2(\Omega)}\,e^{-\sigma^* t} \quad \forall\, t \geq 0. \label{eq52}
\end{align}
Using \eqref{eq30} and Lemma~\ref{lem3}, for $C_\Delta := C_{u_t} + \lambda\,\Gamma$, we obtain
\begin{align}
\|\Delta u(\cdot,t)\|_{L^2(\Omega)} \leq C_\Delta\,\|u_o\|_{H^2(\Omega)}\,e^{-\sigma^* t} \quad \forall\, t \geq 0, \label{eq53}
\end{align}
which, combined with \eqref{eq29}, implies \eqref{eq32}.

\subsection{Proof of \eqref{eq34}-\eqref{eq35}}\label{7}

Differentiating \eqref{eq22} with respect to time, note that we have
\begin{align}
\partial_tv_i(x,t) := \langle M_{\gamma_i}A\dot{U}(t),\mathcal{L}(x)\rangle_{N} \quad  x\in \partial \Omega, ~ t\geq 0. \label{eq37}
\end{align}
Hence,
$$\|\partial_t v_i(\cdot,t)\|_{H^{3/2}(\partial\Omega)} \leq C_{\mathcal{T}}\,\|M_{\gamma_i}A\|\,|\dot{U}(t)|.$$
Since $U$ satisfies $\dot{U} = -SU$ with $S := \sum_{k=1}^{N}\gamma_k B_k A$, we have, according to \eqref{eq28},
$|\dot{U}(t)| \leq \|S\|\,C_1\,|U(0)|\,e^{-\sigma_o t}.$
Furthermore,  
$|U(0)| \leq \|u_o\|_{L^2(\Omega)} \leq \|u_o\|_{H^2(\Omega)}$
holds by Bessel's inequality \cite[Thm.~5.9, pp.~141--142]{brezis}. Now, recalling from \eqref{eq16} that
$\xi_i(\cdot,t) = D_{\gamma_i}(v_i(\cdot,t)),$
we have
$\partial_t \xi_i(\cdot,t) = D_{\gamma_i}(\partial_t v_i(\cdot,t)).$
Hence,
$$\|\partial_t \xi_i(\cdot,t)\|_{H^2(\Omega)} \leq C_{D,\gamma_i}\,\|\partial_t v_i(\cdot,t)\|_{H^{3/2}(\partial\Omega)}.$$
Setting
$$C_\xi := \max_{1 \leq i \leq N}\{C_{D,\gamma_i}\}\,C_{\mathcal{T}}\,
\max_{1 \leq i \leq N}\{\|M_{\gamma_i}A\|\}\,C_1\,\max(1,\|S\|),$$
we obtain \eqref{eq34}. Inequality \eqref{eq35} follows similarly, 
replacing $|\dot{U}(t)| \leq \|S\|\,C_1\,|U(0)|\,e^{-\sigma_o t}$ 
by $|U(t)| \leq C_1\,|U(0)|\,e^{-\sigma_o t}$, which is absorbed 
by the factor $\max(1,\|S\|)$ in $C_\xi$.

\subsection{Proof of \eqref{eq38}}\label{8}

By Parseval's identity \cite[Cor.~5.10, p.~143]{brezis}, it suffices to analyze $\sum_{n=1}^{+\infty}|\dot{w}_n(t)|^2$.

$\bullet$ \textit{Analysis of the first $N$ modes}: differentiating \eqref{eq17} with respect to time, we have 
$\dot{w}_n = \dot{u}_n - \sum_{i=1}^{N}\langle \partial_t \xi_i, \varphi_n\rangle_{L^2(\Omega)}.$
Combining this identity with $\dot{U} = -SU$ and \eqref{eq34}, we obtain 
\begin{align}
\sum_{n=1}^{N}|\dot{w}_n(t)|^2 \leq C_u^2\,\|u_o\|_{H^2(\Omega)}^2\,e^{-2\sigma_o t} \quad \forall\, t \geq 0, \label{eq39}
\end{align}
where $C_u := \sqrt{N}\,(\|S\|\,C_1 + N\,C_\xi)$.

$\bullet$ \textit{Analysis of the remaining modes}: for each $n \geq N+1$, note that we have, according to $\Sigma_w$,
\begin{align}
\dot{w}_n = \mu_n w_n + g_n, \label{eq40}
\end{align}
where $g_n := \langle g, \varphi_n\rangle_{L^2(\Omega)}$ and
$$g := \sum_{i=1}^{N}\gamma_i\xi_i - \sum_{i=1}^{N}\partial_t\xi_i + 2\sum_{k=1}^{N}\mu_k\Bigl(\sum_{i=1}^{N}\langle\xi_i,\varphi_k\rangle_{L^2(\Omega)}\Bigr)\varphi_k.$$
From \eqref{eq34} and \eqref{eq35}, we have
\begin{align}
\|g(\cdot,t)\|_{L^2(\Omega)} \leq C_g\,\|u_o\|_{H^2(\Omega)}\,e^{-\sigma_o t} \quad \forall\, t \geq 0, \label{eq41}
\end{align}
where $C_g := \bigl(\sum_{i=1}^{N}\gamma_i + N \max(1,\|S\|) + 2N\mu_1\bigr)\,C_\xi$. 

Differentiating $g$ with respect to time, we obtain
$$\partial_t g = \sum_{i=1}^{N}\gamma_i\partial_t\xi_i \!-\! \sum_{i=1}^{N}\partial_t^2\xi_i \!+\! 2\!\sum_{k=1}^{N}\!\mu_k\!\Bigl(\sum_{i=1}^{N}\langle\partial_t\xi_i,\varphi_k\rangle_{L^2(\Omega)}\!\Bigr)\!\varphi_k.$$
Since $\partial_t^2\xi_i = D_{\gamma_i}(\partial_t^2 v_i)$, 
$\partial_t^2 v_i(x,t) = \langle M_{\gamma_i}A\ddot{U}(t),\,\mathcal{L}(x)\rangle_N,$
and $\ddot{U} = -S\dot{U}$, as in the proof of Lemma~\ref{lem4}, we obtain 
$$\|\partial_t^2\xi_i(\cdot,t)\|_{L^2(\Omega)} \leq C_\xi'\,\|u_o\|_{H^2(\Omega)}\,e^{-\sigma_o t}$$
for some $C_\xi'>0$ independent of $u_o$. Combining with \eqref{eq34}, 
\begin{align}
\|\partial_t g(\cdot,t)\|_{L^2(\Omega)} \leq C_{\dot{g}}\,\|u_o\|_{H^2(\Omega)}\,e^{-\sigma_o t} \quad \forall\, t \geq 0, \label{eq42}
\end{align}
where
$C_{\dot{g}} := \bigl(\textstyle\sum_{i=1}^{N}\gamma_i + 2N\mu_1\bigr)\,C_\xi + N\,C_\xi'.$
Applying Duhamel's formula to \eqref{eq40} and differentiating with respect to $t$, we obtain
\begin{align*}
\dot{w}_n(t) =&~ \mu_n e^{\mu_n t}w_n(0) + g_n(t) + \int_0^t \mu_n e^{\mu_n(t-s)}g_n(s)\,ds. 
\end{align*}
Using integration by parts, we have
$
\int_0^t \mu_n e^{\mu_n(t-s)}g_n(s)\,ds = -g_n(t) + e^{\mu_n t}g_n(0)  + \int_0^t e^{\mu_n(t-s)}\dot{g}_n(s)\,ds. 
$
Hence, we obtain
\begin{align*}
\dot{w}_n(t) =&~ \mu_n e^{\mu_n t}w_n(0) + e^{\mu_n t}g_n(0) + \int_0^t e^{\mu_n(t-s)}\dot{g}_n(s)\,ds. 
\end{align*}
We now bound each of the three terms at the right-hand side of the latter equation. Since $w(\cdot,0) \in D(\mathcal{A})$, Parseval's identity yields
\begin{align}
\sum_{n \geq N+1}\mu_n^2\,e^{2\mu_n t}|w_n(0)|^2 \leq C_0^2\,\|u_o\|_{H^2(\Omega)}^2\,e^{-2|\mu_{N+1}| t}, \label{eq46}
\end{align}
where $C_0 := (\sqrt{n}+\lambda)(1 + N\,C_\xi)$, with $n$ the dimension of $\Omega$. Similarly, by Parseval's identity and \eqref{eq41} at $t=0$, we have
\begin{align}
\sum_{n \geq N+1}e^{2\mu_n t}|g_n(0)|^2 \leq C_g^2\,\|u_o\|_{H^2(\Omega)}^2\,e^{-2|\mu_{N+1}| t}. \label{eq47}
\end{align}
For the third term, we use Minkowski's integral inequality \cite[p.~707]{Evans2010} to obtain 
\begin{align*}
&\left(\sum_{n \geq N+1}\left|\int_0^t e^{\mu_n(t-s)}\dot{g}_n(s)\,ds\right|^2\right)^{1/2} \nonumber\\ &\qquad\leq \int_0^t e^{-|\mu_{N+1}|(t-s)}\|\dot{g}(\cdot,s)\|_{L^2(\Omega)}\,ds. 
\\ &\qquad\leq C_{\dot{g}}\|u_o\|_{H^2(\Omega)}e^{-|\mu_{N+1}|t}\int_{0}^{t}e^{(-\sigma_o+|\mu_{N+1}|)s}ds,
\end{align*}
where the latter inequality substitutes \eqref{eq42}.

If $\sigma_o = |\mu_{N+1}|$, then 
\begin{align*}
&e^{-|\mu_{N+1}|t}\int_{0}^{t}e^{(-\sigma_o+|\mu_{N+1}|)s}ds = te^{-|\mu_{N+1}|t}\\
&\hspace{2.6cm}\leq \frac{e^{-\varepsilon |\mu_{N+1}|t}}{(1-\varepsilon)e|\mu_{N+1}|}
\end{align*}
for any $\varepsilon \in (0,1)$. If $\sigma_o > |\mu_{N+1}|$, then 
\begin{align*}
e^{-|\mu_{N+1}|t}\int_{0}^{t}e^{(-\sigma_o+|\mu_{N+1}|)s}ds \leq \frac{e^{-|\mu_{N+1}|t}}{\left(\sigma_o -|\mu_{N+1}| \right)}.
\end{align*}
If $\sigma_o < |\mu_{N+1}|$, then 
$e^{-|\mu_{N+1}|t}\int_{0}^{t}e^{(-\sigma_o+|\mu_{N+1}|)s}ds \leq e^{-\sigma_o t}/(|\mu_{N+1}|-\sigma_o)$. 
Hence, there exist $\bar{C}_{\dot{g}},\sigma_1> 0$ such that 
\begin{align}
\hspace{-0.4cm} \sum_{n \geq N+1}\left|\int_0^t e^{\mu_n(t-s)}\dot{g}_n(s)\,ds\right|^2 \leq \bar{C}_{\dot{g}}^2\|u_o\|_{H^2(\Omega)}^2e^{-2\sigma_1 t}. \label{eq50}
\end{align}
Combining \eqref{eq46}, \eqref{eq47}, and \eqref{eq50}, we conclude that there exist $C_s, \sigma^*>0$ such that
\begin{align}
\sum_{n \geq N+1}|\dot{w}_n(t)|^2 \leq C_s^2\,\|u_o\|_{H^2(\Omega)}^2\,e^{-2\sigma^* t} \quad \forall\, t \geq 0. \label{eq51}
\end{align}
Combining \eqref{eq39} and \eqref{eq51} and using Parseval's identity yields \eqref{eq38} with $C_w := \sqrt{C_u^2 + C_s^2}$, which concludes the proof.

\section{Simulation Results}\label{sec6}

We illustrate our results with numerical simulations performed in Matlab, for $\Omega = \{x\in \mathbb{R}^n: |x|< 2\}$ with $n=2$ and $n=3$. We set $\lambda = 6.6141$, which yields $N=5$ nonnegative eigenvalues for $n=2$ and $N=4$ for $n=3$, and $u_o(x):=(4-|x|^2)p_3(x)$, where $p_3$ is a third-degree polynomial with randomly generated coefficients in $[-1,1]$. The PDE is discretized by projecting it onto the first $300$ eigenfunctions of $\mathcal{A}$, and the resulting ODE is solved over $[0,4]$ with the solver \texttt{ode15s} and a time step $0.05$. The $H^2$ norm is approximated by $\big(\sum_{n}(1 + \mu_n^2)u_n(t)^2\big)^{1/2}$, and the max norm on a uniform grid of $50^n$ points.

\begin{figure}
\centering
\includegraphics[width=0.235\textwidth]{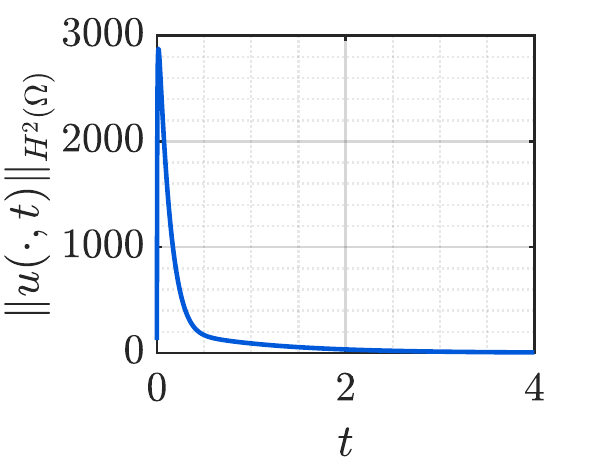}\hfill
\includegraphics[width=0.235\textwidth]{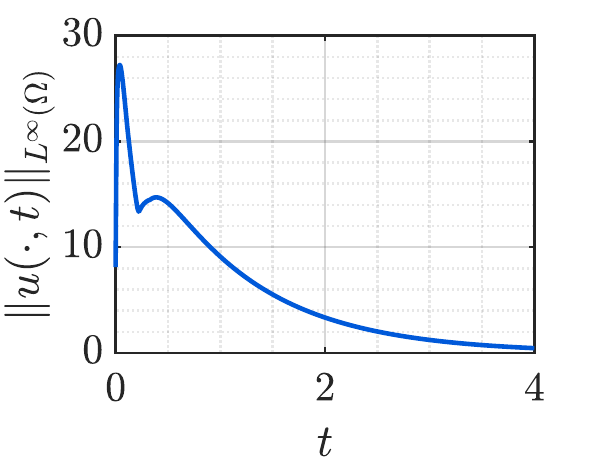}
\vspace{-0.6cm}
\caption{
The $2$-D case.}
\label{fig2}
\end{figure}

In the $2$-D case we found $\mu_1 \approx 5.17$, $\mu_2 = \mu_3 \approx 2.94$, $\mu_4 = \mu_5 \approx 0.02$, and we select $\gamma_i=\mu_1+i$. Figure~\ref{fig2} plots $t\mapsto \|u(\cdot,t)\|_{H^2(\Omega)}$ and $t\mapsto \|u(\cdot,t)\|_{L^{\infty}(\Omega)}$, which converge exponentially to zero, so the state converges to zero at every spatial location. When $v=0$, both norms diverge.

\begin{figure}
\centering
\includegraphics[width=0.235\textwidth]{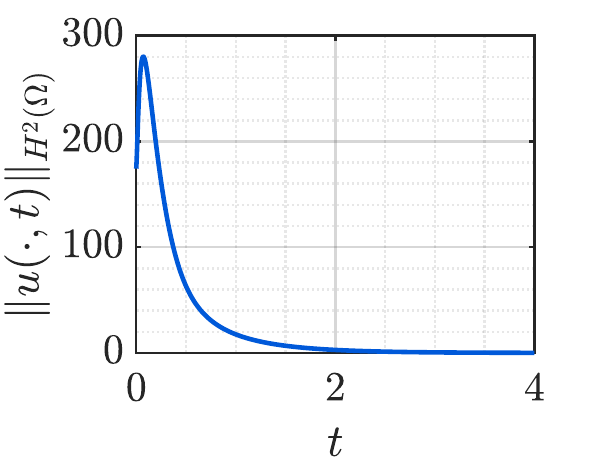}\hfill
\includegraphics[width=0.235\textwidth]{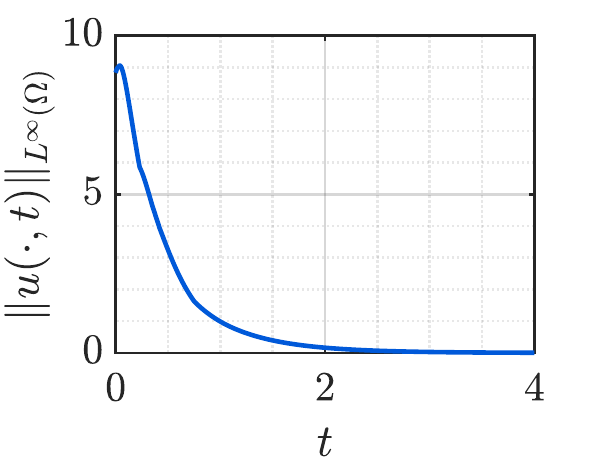}
\vspace{-0.6cm}
\caption{
The $3$-D case.}
\label{fig5}

\vspace{-0.4cm}
\end{figure}

Next, in $3$-D, we found $\mu_1 \approx 4.147$ and $\mu_2 = \mu_3 = \mu_4 \approx 1.566$, and we set $\gamma_i=\mu_1+i$. Figure~\ref{fig5} shows exponential decay of both norms to zero. 

\section{Conclusion and Perspectives} \label{sec7}
We proved $H^2$ exponential stability of the origin for linear parabolic equations in $2$-D and $3$-D, using the boundary controller in  \cite{Munteanu2017IJC}. Our result implies boundedness and pointwise convergence of the state's max norm to zero. Our approach uses the 
previously-established $H^1$ exponential stability of the origin, transforming the problem into the analysis of the $L^2$ norm of the state's Laplacian, which, in turn, reduces to analyzing the $L^2$ norm of the state's time derivative.

Several research perspectives can be envisioned. For instance, in the 
output-feedback setting, the controller is coupled with an observer, and the analysis of the time derivative of the vector of unstable modes, needs to be revisited to account for the estimation error. Moreover, the boundary is assumed to be of class $\mathcal{C}^3$, and the case of nonsmooth domains remains open.

\appendix

\subsection{Proof of Lemma \ref{lem_D}}\label{p_lem_D}
Item (i) is in \cite[Sec.~3.1]{lhachemi2025automatica}. Let us prove item (ii). Linearity is immediate from the linearity of $\Sigma_H$ in $(D,f)$. To prove continuity, we show that the graph of $D_\gamma:H^{3/2}(\partial\Omega)\to H^2(\Omega)$ is closed and invoke the closed graph theorem \cite[Thm.~2.9, p.~37]{brezis}. Let $(f_m)_{m\geq 1}\subset H^{3/2}(\partial\Omega)$ and set $D_m:=D_\gamma(f_m)\in H^2(\Omega)$. Assume $f_m\to f$ in $H^{3/2}(\partial\Omega)$ as $m\to +\infty$, and $D_m\to D$ in $H^2(\Omega)$ as $m\to +\infty$. Since $\Omega$ is a $\mathcal{C}^3$ domain, the Dirichlet trace operator $\omega_0:H^2(\Omega)\to H^{3/2}(\partial\Omega)$, $\omega_0(f) := f|_{\partial \Omega}$, is continuous. Hence $\omega_0(D)=\lim_{m\to\infty}\omega_0(D_m) =\lim_{m\to\infty} f_m=f$. Next, the operator $\bar{\mathcal{A}}: H^2(\Omega)\to L^2(\Omega)$ is continuous, and the finite-dimensional operator $D\mapsto \sum_{i=1}^{N} \langle D,\varphi_i\rangle_{L^2(\Omega)}\,\varphi_i$ is bounded on $L^2(\Omega)$. Therefore, the left-hand side of the first equation in $\Sigma_H$ defines a continuous map from $H^2(\Omega)$ to $L^2(\Omega)$. Passing to the limit as $m\to +\infty$, we conclude that $D$ solves $\Sigma_H$. By uniqueness, $D=D_\gamma(f)$. Thus the graph of $D_\gamma$ is closed, and the closed graph theorem implies that $D_\gamma$ is bounded, hence continuous.

    Finally, we prove item (iii). Fix $t \geq 0$ and $h \neq 0$. By the linearity of $D_\gamma$,
$$
\frac{D_\gamma(v(\cdot,t+h)) - D_\gamma(v(\cdot,t))}{h} 
= D_\gamma\!\left(\frac{v(\cdot,t+h) - v(\cdot,t)}{h}\right).
$$
Since $v \in \mathcal{C}^1([0,+\infty); H^\frac{3}{2}(\partial\Omega))$, the argument on the right-hand side converges to $\partial_t v(\cdot,t)$ in $H^{\frac{3}{2}}(\partial\Omega)$ as $h \to 0$. By the continuity of $D_\gamma$, the right-hand side converges to $D_\gamma(\partial_t v(\cdot,t))$ in $H^2(\Omega)$. Hence $\partial_t D_\gamma(v(\cdot,t)) = D_\gamma(\partial_t v(\cdot,t))$. Continuity of $t \mapsto D_\gamma(\partial_t v(\cdot,t))$ in $H^2(\Omega)$ follows from the continuity of $t\mapsto \partial_t v(\cdot,t)$ and of $D_\gamma$.


\begin{thebibliography}{99}
\footnotesize
\bibitem{moh2} M.~A.~Ouchdiri, S.~Benjelloun, A.~Saoud, and I.~Otero-Muras, ``Turing patterns in a morphogenetic model with single regulatory function,'' \emph{Mathematical Biosciences}, vol.~389, Art.~no.~109536, 2025.
\bibitem{moh} M.~A.~Ouchdiri, H.~Faquir, S.~Benjelloun, M.~Maghenem, I.~Otero-Muras, and A.~Saoud, ``An optimal-control framework for reaction diffusion systems with application to synthetic developmental biology,'' in \textit{Proc. IEEE 64th Conf. Decision and Control}, pp.~1925--1930, 2025.
\bibitem{LasieckaTriggiani2000} I.~Lasiecka and R.~Triggiani, \emph{Control Theory for Partial Differential Equations, Vol.~I: Abstract Parabolic Systems}. Cambridge University Press, 2000.
\bibitem{VazquezKrstic2016} R.~Vazquez and M.~Krsti\'c, ``Explicit output-feedback boundary control of reaction-diffusion PDEs on arbitrary-dimensional balls,'' \emph{ESAIM: COCV}, vol.~22, no.~4, pp.~1078--1096, 2016.
\bibitem{LiuXie2020} X.~Liu and C.~Xie, ``Boundary control of reaction--diffusion equations on higher-dimensional symmetric domains,'' \emph{Automatica}, vol.~114, Art.~no.~108832, 2020.
\bibitem{Meurer2013} T.~Meurer, \emph{Control of Higher-Dimensional PDEs: Flatness and Backstepping Designs}. Springer, 2013.
\bibitem{higher_1}
M. C. Belhadjoudja, ''Backstepping Control of Parabolic PDEs on Directionally-Convex Domains,'' \emph{Under review}. Preprint Available at SSRN 6447745.

\bibitem{higher_2}
R. Vazquez, ''Output-Feedback Boundary Control of Reaction--Diffusion PDEs on Arbitrary Lipschitz Domains: A Target-Domain Approach,''
\emph{arXiv preprint:} arXiv:2606.16345.

\bibitem{Balas1978} M.~J.~Balas, ``Active control of flexible systems,'' \emph{J. Optim. Theory Appl.}, vol.~25, no.~3, pp.~415--436, 1978.
\bibitem{KatzFridman2020} R.~Katz and E.~Fridman, ``Constructive method for finite-dimensional observer-based control of 1-D parabolic PDEs,'' \emph{Automatica}, vol.~122, Art.~no.~109285, 2020.
\bibitem{KatzFridman2021} R.~Katz and E.~Fridman, ``Delayed finite-dimensional observer-based control of 1-D parabolic PDEs,'' \emph{Automatica}, vol.~123, Art.~no.~109364, 2021.
\bibitem{LhachemiPrieur2022} H.~Lhachemi and C.~Prieur, ``Finite-dimensional observer-based boundary stabilization of reaction--diffusion equations with either a Dirichlet or Neumann boundary measurement,'' \emph{Automatica}, vol.~135, Art.~no.~109955, 2022.
\bibitem{Barbu2013TAC} V.~Barbu, ``Boundary stabilization of equilibrium solutions to parabolic equations,'' \emph{IEEE Trans. Autom. Control}, vol.~58, no.~9, pp.~2416--2420, 2013.
\bibitem{Munteanu2017IJC} I.~Munteanu, ``Stabilisation of parabolic semilinear equations,'' \emph{Int. J. Control}, vol.~90, no.~5, pp.~1063--1076, 2017.
\bibitem{lhachemi2025automatica} H.~Lhachemi, I.~Munteanu, and C.~Prieur, ``Boundary output feedback stabilization for 2-D and 3-D parabolic equations,'' \emph{Automatica}, vol.~176, Art.~no.~112259, 2025.
\bibitem{lhachemi2024stabilization} H.~Lhachemi and C.~Prieur, ``Stabilization of a reaction-diffusion equation in $H^2$-norm with application to saturated Neumann measurement,'' in \textit{Proc. IEEE 63rd Conf. Decision and Control}, pp.~1187--1192, 2024.
\bibitem{LhachemiPrieur2026H2} H.~Lhachemi and C.~Prieur, ``Output feedback stabilization of reaction--diffusion equations in $H^2$-norm,'' \emph{IEEE Trans. Autom. Control}, 2026, doi:10.1109/TAC.2026.3688198.
\bibitem{Evans2010} L.~C.~Evans, \emph{Partial Differential Equations}, 2nd ed. American Mathematical Society, 2010.
\bibitem{Grisvard} P.~Grisvard, \emph{Elliptic Problems in Nonsmooth Domains}. Pitman, Boston, 1985.
\bibitem{brezis} H.~Brezis, \emph{Functional Analysis, Sobolev Spaces and Partial Differential Equations}. Springer, New York, 2011.
\end{thebibliography}
\end{document}